\documentclass[11pt]{amsart}

\newcommand{\BF}{{\mathbf F}}
\newcommand{\BQ}{{\mathbf Q}}
\newcommand{\BR}{{\mathbf R}}
\newcommand{\BZ}{{\mathbf Z}}

\newcommand{\pid}{{\mathfrak p}}

\newcommand{\eps}{\epsilon}
\newcommand\ra{\rightarrow}

\DeclareMathOperator{\Tr}{Tr}

\newtheorem{theorem}{Theorem}
\newtheorem{lemma}[theorem]{Lemma}

\newtheorem*{assumption}{Assumption}

\theoremstyle{definition}

\theoremstyle{remark}

\newtheorem{acks}{Acknowledgments}	
\newtheorem{conventions}{Conventions}	

\begin{document}

\title{Higher-order Carmichael numbers}
\author{Everett W.\ Howe}
\address{Center for Communications Research, 
         4320 Westerra Court, 
         San Diego, CA 92121-1967, USA.}
\email{however@alumni.caltech.edu}
\urladdr{http://alumni.caltech.edu/\~{}however/}

\date{2 December 1998}


\begin{abstract}
We define a Carmichael number of order $m$ to be a composite integer $n$ such that 
$n$th-power raising defines an endomorphism of every $\BZ/n\BZ$-algebra that can be
generated as a $\BZ/n\BZ$-module by $m$ elements.  We give a simple criterion to determine whether
a number is a Carmichael number of order~$m$, and we give a heuristic argument
(based on an argument of Erd\H{o}s for the usual Carmichael numbers)
that indicates that for every $m$ there should be infinitely many Carmichael numbers
of order~$m$.  The argument suggests a method for finding examples of higher-order
Carmichael numbers; we use the method to provide examples of Carmichael numbers of
order~$2$.
\end{abstract}

\maketitle

\section{Introduction}
\label{S-intro}

A Carmichael number is defined to be
a positive composite integer $n$ that is a Fermat pseudoprime to every base;
that is, a composite $n$ is a Carmichael number if $a^n\equiv a\bmod n$ for every integer $a$.
Clearly one can generalize the idea of a Carmichael number by 
allowing the pseudoprimality test in the definition to vary over some larger class of tests
(perhaps including some of those found in \cite{adams}, \cite{adams-shanks}, \cite{baillie-wagstaff},
\cite{diporto-filipponi-montolivo}, \cite{grantham}, \cite{gurak}, \cite{joo},
\cite{lidl-muller-generalizations}, \cite{marko}, \cite{szekeres}),
and indeed such generalizations have been considered (see for example
\cite{diporto-filipponi}, \cite{grantham}, \cite{kowol}, 
\cite{lidl-muller-note}, \cite{lidl-muller-primality}, \cite{lidl-muller-oswald}, 
\cite{marko}, \cite{s-muller}, \cite{muller-oswald}, \cite{williams}).
But there is also a natural algebraic way of generalizing the concept of 
a Carmichael number that makes no mention of pseudoprimality.
To motivate the definition we note that 
(1) an integer $n>1$ is prime if and only if $n$th-power
    raising is an endomorphism of every $\BZ/n\BZ$-algebra, and 
(2) a positive composite integer $n$ is a Carmichael number if 
    and only if $n$th-power raising is an endomorphism of $\BZ/n\BZ$.
So if $m$ is a positive integer, we define a
{\it Carmichael number of order $m$\/} to be a positive composite integer $n$
such that the function $x\mapsto x^n$ defines an endomorphism of every
$\BZ/n\BZ$-algebra that can be generated as a $\BZ/n\BZ$-module by $m$ elements.

Although our definition does not explicitly mention pseudoprimality, a
Carmichael number $n$ of order $m$ will pass many reasonable pseudoprimality tests.
For example, if $\alpha$ is an algebraic integer of degree $d$ with $d\le m$,
then we have $\Tr_{\BQ(\alpha)/\BQ} (\alpha^n) \equiv \Tr_{\BQ(\alpha)/\BQ} (\alpha) \bmod n$,
so $n$ will pass a Dickson-like pseudoprimality test based on 
the recurrence sequence of order $d$ consisting of the traces of the powers of $\alpha$.
Also, $n$ will pass the ``Frobenius step'' of the Frobenius pseudoprime
test of Grantham~\cite{grantham} with respect to every polynomial of degree at most $m$.

We will prove the following theorem, 
which provides a characterization of the Carmichael numbers of order $m$
that generalizes Korselt's criterion~\cite{korselt} for the usual Carmichael numbers:

\begin{theorem}
\label{characterization}
Let $m$ and $n$ be positive integers with $n$ composite.
The following statements are equivalent:
\begin{enumerate}
\item[(a)] $n$ is a Carmichael number of order $m$\textup{;}
\item[(b)] the function $x\mapsto x^n$ defines an endomorphism of every 
finite \'etale $\BZ/n\BZ$-algebra that can be generated as a $\BZ/n\BZ$-module by $m$ elements\textup{;}
\item[(c)] the following two conditions hold\textup{:}
\begin{enumerate}
\item[(i)] $n$ is squarefree\textup{;}
\item[(ii)] for every prime divisor $p$ of $n$ and for every integer $r$ with $1\le r\le m$,
there is an integer $i\ge 0$ such that $n\equiv p^i \bmod (p^r-1)$.
\end{enumerate}
\end{enumerate}
\end{theorem}

(For the benefit of those readers unfamiliar with finite \'etale $R$-algebras, we present
a definition equivalent to the usual one (found for example in Section~I.3 of~\cite{milne})
that is applicable when $R$ is a finite product of local rings.
First suppose that $R$ is itself a local ring --- that is, a ring with a unique maximal ideal.  
Then an $R$-algebra $S$ is {\it finite \'etale\/} if it is
free of finite rank as an $R$-module and if for some (or equivalently, every)
$R$-module basis $\{e_1,\ldots,e_n\}$ of $S$, the determinant of the $n$-by-$n$ matrix
$[\Tr_{S/R}(e_i e_j)]$ is a unit of $R$; here $\Tr_{S/R}$ is the
trace map from $S$ to $R$.
Now suppose $R = R_1\times\cdots\times R_m$, where the $R_i$ are local.
Then an $R$-algebra $S$ is {\it finite \'etale\/} if it is of the form
$S = S_1\times \cdots\times S_m$, where each $S_i$ is a finite \'etale
$R_i$-algebra.  (Note that the zero ring is a finite \'etale $R_i$-algebra,
so some of the $S_i$ may be zero.) 
Since every finite ring is a finite product of local rings,
our definition can be used when $R$ is finite.
We see, for example, that if $n$ is a squarefree integer then
a finite \'etale $\BZ/n\BZ$-algebra is simply a finite product of finite fields,
each of whose characteristics divides $n$.)

Theorem~\ref{characterization} allows us to formulate
a heuristic argument (based on an argument of Erd\H{o}s~\cite{erdos} for the
usual Carmichael numbers, 
and similar to an argument of Pomerance~\cite{pomerance} for the
Baillie-PSW pseudoprimes) that indicates that for every $m$ there should be
infinitely many Carmichael numbers of order~$m$.  
The heuristics suggest a method of searching for higher-order Carmichael numbers;
we implement this method for the case $m=2$ and find many examples, some of which we present below.
In fact, the numbers $n$ produced by our argument have the property that 
$n$th-power raising is the {\it identity\/} on every finite \'etale $\BZ/n\BZ$-algebra
that can be generated as a module by $m$ elements.
We call such $n$ {\it rigid\/} Carmichael numbers of order $m$, 
and in Section~\ref{S-examples} we show by example that not all higher-order
Carmichael numbers are rigid.

We would like to replace the heuristic arguments of this paper with actual proofs,
but that seems to be difficult; we have been unable to adapt the argument of 
Alford, Granville, and Pomerance~\cite{alford-granville-pomerance}
for the infinitude of the usual Carmichael numbers to the case of
higher-order Carmichael numbers.
However, in a recent paper~\cite{hsu}, Hsu proves that there are infinitely many
``Carmichael polynomials'', which are Drinfeld module analogues of Carmichael
numbers and higher-order Carmichael numbers.

We know of only one example of a higher-order Carmichael number other
than the ones produced by the computations described in this paper:
one finds the number $17 \cdot 31 \cdot 41 \cdot 43 \cdot 89 \cdot 97 \cdot 167 \cdot 331$,
which is a rigid Carmichael number of order~$2$, on the list of 
the Carmichael numbers less than~$10^{16}$ that was computed by Richard Pinch
(see~\cite{pinch}, \cite{pinch-ftp}).

\begin{acks} 
The author thanks Dan Gordon, Jon Grantham, Andrew Granville, Hendrik Lenstra, and Carl Pomerance
for reading and commenting on various versions of this note.
The author is especially grateful to Lenstra for suggesting 
Lemma~\ref{smallprimes} and its proof, and for suggesting various ways of
defining ``finite \'etale'' without using much algebra.
\end{acks}

\begin{conventions}
We subscribe to the conventions that rings have 
identity elements and that ring homomorphisms $R\ra S$ take
the identity of $R$ to the identity of $S$.
\end{conventions}

\section{Proof of the Theorem}
\label{S-proof}

The implication (a) $\Rightarrow$ (b) is trivial.

Suppose that condition (b) holds.  The ring $\BZ/n\BZ$ is a finite \'etale
algebra over itself and is generated by a single element as a module over itself, 
so $x\mapsto x^n$ must be an endomorphism of this ring.
The only endomorphism of $\BZ/n\BZ$ is the identity,
so we have $x = x^n$ for all $x$ in $\BZ/n\BZ$.  But if $n$ were divisible by the square of
a prime $p$ we would have $p^n \not\equiv p \bmod n$, a contradiction.
Thus $n$ is squarefree.

Let $p$ be a prime divisor of $n$ and let $r$ be an integer with
$1\le r\le m$. Let $F$ be the finite field with 
$p^r$ elements.  The field $F$ is a finite \'etale $\BZ/p\BZ$-algebra, and so is also a
finite \'etale $\BZ/n\BZ$-algebra via the projection $\BZ/n\BZ\ra \BZ/p\BZ$.  It is clear that $F$
can be generated as a $\BZ/n\BZ$-module by $m$ elements, so $n$th-power raising
is an automorphism of $F$.  Every automorphism of $F$ is of the form $x\mapsto x^{p^i}$
for some $i$, so there is an integer $i$ such that $x^n = x^{p^i}$ for every $x\in F$.
Since the multiplicative group of $F$ is cyclic of order $p^r-1$, we see that
$n \equiv p^i \bmod (p^r-1)$.  This proves the implication (b) $\Rightarrow$ (c).

Now suppose that condition (c) holds.  First we prove the following statement:

\begin{lemma}
\label{smallprimes}
If $r$ is an integer with $1 \le r \le m$ then $\binom{n}{r} \equiv 0 \bmod n$.
\end{lemma}

\begin{proof}
Note that the statement we are to prove
is equivalent to the statement that all prime divisors of $n$ are greater than~$m$.
Suppose, to obtain a contradiction, that $n$ had a prime divisor $q$ with $q\le m$.
Since $n$ is assumed to be composite and squarefree, $n$ must have another prime
divisor $p\neq q$.  If we apply statement (c)(ii) of the theorem with $r = q - 1$,
we find that $n\equiv p^i \bmod (p^{q-1}-1)$, and since $q$ divides $p^{q-1}-1$
it follows that $n \equiv p^i \bmod q$.
But $q \mid  n$, so we find that $q\mid p^i$, a contradiction.
\end{proof}

Now suppose $R$ is a $\BZ/n\BZ$-algebra that can be generated as a module
by $m$ elements.  Then $R$ is a finite ring, and so is a product of
finite local rings $R_i$, each of which is a $\BZ/n\BZ$-algebra that
can be generated as a $\BZ/n\BZ$-module by $m$ elements.  
If $n$th-power raising is an endomorphism of each $R_i$, then it is an endomorphism
of $R$ as well, so it suffices to consider the case where $R$ is local.  
Since $n$ is squarefree, there is a prime divisor $p$ of $n$ such that
$pR = 0$, so that $R$ is an $\BF_p$-algebra.  
Let $\pid$ be the maximal ideal of $R$ and let $k = R/\pid$.
Since $R$ can be generated by $m$ elements as an $\BF_p$-module,
we see that $[k:\BF_p]\le m$ and that $\pid^m = 0$.
Since $k$ is separable over $\BF_p$,
Hensel's lemma shows that there is a homomorphism $k\ra R$
compatible with the reduction map $R\ra k$; we view $k$ as a subring of $R$ via this map.
We find that every element of $R$ may be written in a unique way as a sum
$a + z$ where $a\in k$ and $z\in\pid$.

If $a\in k$ and $z\in\pid$, then we have
$$(a+z)^n = \sum_{r=0}^n \binom{n}{r} a^{n-r} z^r = a^n$$
where the second equality is 
obtained from the facts that  $z^r = 0 $ when $r\ge m$ and $\binom{n}{r} = 0$ in $R$ when $1\le r\le m$.
But since $n\equiv p^i \bmod (p^{[k:\BF_p]} - 1)$ we see that
$(a+z)^n = a^{p^i}$,
so $n$th-power raising on $R$ is simply the reduction map to $k$ followed by the
automorphism $x\mapsto x^{p^i}$ followed by the lifting map $k\ra R$.
In particular, $n$th-power raising is a homomorphism.  Thus, $n$ is a Carmichael
number of order $m$.  This shows that (c) implies (a), and completes the proof of the theorem.

\section{A construction and heuristics}
\label{S-heuristics}

Let $m>0$ be given. In this section we will give a construction that associates to
every positive integer $L$ a (possibly empty) set $C(m,L)$ of Carmichael numbers of
order $m$.  We will also give a heuristic argument that indicates that one should be able to find
values of $L$ that will make $\#C(m,L)$ as large as one pleases.
The construction and argument generalize those of Erd\H{o}s~\cite{erdos} for the usual Carmichael numbers;
Pomerance uses a similar argument in~\cite{pomerance}
to show that there should be infinitely many Baillie-PSW pseudoprimes.

First, the construction.  Let $P(m,L)$ be the set of prime numbers $p$ that do not
divide $L$ and that have the property that for every positive integer $r\le m$,
the integer $p^r-1$ divides $L$.  Let $C(m,L)$ be the set of squarefree integers $n>1$
that are congruent to $1$ modulo $L$ and whose prime divisors all lie in $P(m,L)$.
We claim that the elements of $C(m,L)$ are Carmichael numbers of order $m$. For suppose
$n\in C(m,L)$, suppose $r$ is an integer with $1\le r\le m$, and suppose $p$ is a prime
divisor of $n$.  Then $p^r-1$ divides $L$, and $L$ divides $n-1$, so 
$n\equiv p^0\bmod (p^r-1)$.  By Theorem~\ref{characterization}, the 
integer $n$ is a Carmichael number of order $m$.

Our heuristic argument for the existence of $L$ for which $\#C(m,L)$ is large
depends on the following assumption (in addition to the usual assumptions and approximations
made in such arguments):

\begin{assumption}
Suppose $f$ is an element of $\BZ[x]$ and $\eps$ is a positive real.
Then there is a positive integrable function $s$ from $[1,1+\eps]$ to $\BR$
such that for $y$ sufficiently large and for every $u\in[1,1+\eps]$
there are at least $y^u s(u)$ integers $x$ in $[1,y^u]$ such that $f(x)$ is $y$-smooth.
\end{assumption}

Let $\eps>0$ be fixed for the remainder of the argument.
Let $y$ be given, let $f$ be 
the least common multiple of the polynomials $x^r-1$ for $1\le r\le m$,
and let $L$ be the least common multiple of the
prime powers $p^e$ such that $p<y$ and $p^e < y^{m(1+\eps)}.$ 
We will argue that one should expect $\log \#C(m,L) \gg y^{1+\eps/2}$.

Let us estimate the cardinality of the set $S(y,\eps)$ of primes $q$
between $y$ and $y^{1+\eps}$ such that $f(q)$ is $y$-smooth.
By our assumption above, there is a positive integrable
function $s$ such that the probability that a randomly-chosen
integer $x$ less than $y^u$ has $f(x)$ being $y$-smooth is at least $s(u)$.
Thus we expect that the probability that a randomly-chosen
integer $x$ near $y^u$ has $f(x)$ being $y$-smooth is also
at least $s(u)$, so it seems reasonable to approximate a lower bound for $\#S(y,\eps)$ by
$$\int_y^{y^{1+\eps}} s(\log x/\log y) \frac{1}{\log x} \ dx.$$
By setting $u=\log x/\log y$ we convert this last integral to
$$\int_1^{1+\eps} \frac{s(u)}{u} y^u \ du.$$
Thus we expect that
$$\#S(y,\eps) > \int_{1+\eps/2}^{1+\eps} \frac{s(u)}{u} y^u \ du
            > y^{1+\eps/2} \int_{1+\eps/2}^{1+\eps} \frac{s(u)}{u}\ du.$$
Let $c_\eps$ denote the rightmost integral, which is nonzero because $s$ is positive.

Suppose $q$ is an element of $S(y,\eps)$ and let $r$ be an integer
with $1\le r\le m$.  Since $f(q)$ is $y$-smooth, we see that
all of the prime factors of $q^r-1$ are less than $y$. 
Suppose $p$ is a prime divisor of $q^r-1$ and suppose $p^e$ is the largest
power of $p$ that divides $q^r-1$.  Then certainly $p^e\le q^r-1 < q^m \le y^{m(1+\eps)}$,
so $p^e$ divides $L$.  It follows that $q^r-1$ divides $L$.  Thus $S(y,\eps)$ is
contained in $P(m,L)$, and $\#P(m,L) > c_\eps y^{1+\eps/2}$.

Consider the map from the power set of $P(m,L)$ to $(\BZ/L\BZ)^*$ defined by
sending a subset of $P(m,L)$ to the residue modulo $L$ of the product of its
elements.  It seems reasonable to assume that the elements of $(\BZ/L\BZ)^*$
will each have roughly the same number of preimages in the power set of $P(m,L)$,
so we expect that there should be roughly $2^{\#P(m,L)}/\varphi(L)$ 
subsets $X$ of $P(m,L)$ such that the product the elements of $X$ 
is $1$ modulo $L$.  In other words, we expect 
$$\log \#C(m,L) \approx \#P(m,L)\log 2 - \log\varphi(L).$$
Now, $\log L$ should be roughly $m y (1+\eps)$, so $\log\varphi(L)$
should be less than that same amount.  It follows that
we should have $\log \#C(m,L) \gg y^{1+\eps/2}$,
and so we expect to be able to find integers $L$ for which $\#C(m,L)$ is a large
as we like.

\section{Constructing Carmichael numbers of order $2$}
\label{S-constructing}

The argument given in Section~\ref{S-heuristics} suggests a method for finding
Carmichael numbers of order $m$:  Find a value of $L$ for which 
$\#P(m,L)\log 2 - \log\varphi(L)$ is large, and then search for subsets of $P(m,L)$
the products of whose elements are $1$ modulo $L$.
Only about $1$ out of every $\varphi(L)$ subsets of $P(m,L)$ will have the desired
property, so if $L$ is too large we will have trouble finding such subsets.
If $m$ is greater than $2$, we must take $L$ to be extremely large in order for 
our heuristics to predict that $C(m,L)$ is nonempty, so 
examples of Carmichael numbers of order~$3$ or more seem to be out of reach
for the moment.
However, as we will show in this section, 
it is possible to use the above method 
to find Carmichael numbers of order~$2$.  

Let us define the {\it fecundity\/} of a number $L$ to be $F(L) = \#P(2,L) - (\log\varphi(L))/\log 2$,
so that we expect $C(2,L)$ to contain about $2^{F(L)}$ elements.  When $L$ does not have
too many divisors, 
one can compute the set $P(2,L)$ na\"{\i}vely by listing the divisors $d$ of $L$ and
searching for those $d$ such that $d+1$ is the square of a prime.
We computed $F(L)$ by this method for many $L$ built up of primes less than or equal to
$37$, and we found several $L$ with positive fecundity.  For example, let 
$$L_1 = 2^7 \cdot 3^3 \cdot 5^2  \cdot 7 \cdot 11 \cdot 13 \cdot 17 \cdot 19 \cdot 29$$
and 
$$L_2 = 2^7 \cdot 3^3 \cdot 5^2 \cdot 7 \cdot 11 \cdot 13 \cdot 17 \cdot 19 \cdot 29 \cdot 31.$$
Then $\#P(2,L_1) = 45$ and $\#P(2,L_2) = 58$, so that
$F(L_1)\approx 8.039$ and $F(L_2)\approx 16.132$.

We used a ``meet-in-the-middle'' approach to find the
elements of $C(2,L_1)$, using the mathematics package MAGMA on one $195$-MHz 
MIPS R10000 IP27 processor of a Silicon Graphics Origin 2000 computer.
In particular, we divided the set $P(2,L_1)$ into three disjoint subsets $S_1$, $S_2$, and $S_3$
with $\#S_1 = \#S_2 = 19$ and $\#S_3 = 7$, and for each $i=1,2,3$ we let $m_i$ be the product 
of the primes in $S_i$.
We calculated the set $X$ of multiplicative inverses
of the residues (modulo $L_1$) of the  $2^{19}$ divisors of $m_1$
and the set $Y$ of the residues (modulo $L_1$) of the $2^{19}$ divisors of $m_2$.
For every one of the $2^7$ divisors $d$  of $m_3$
we calculated the set $Y_d = \{dy : y \in Y\}$.  For every element $x$ in the intersection $X\cap Y_d$,
we found all divisors $e$ of $m_1$ such that $e\equiv x^{-1}\bmod L_1$ and all divisors
$f$ of $m_2$ such that $df \equiv x\bmod L_1$.  
For each such triple $(d,e,f)$ the product $def$ is congruent to $1$ modulo $L_1$,
and so is an element of $C(2,L_1)$ (unless $d = e = f = 1$).
We found that $\#C(2,L_1) = 246$, whereas 
our heuristic argument suggested that there would be approximately $2^{F(L_1)} \approx 263$
elements in this set.
The two elements of $C(2,L_1)$ with the smallest number of prime divisors are 
$$31\cdot37\cdot101\cdot103\cdot109\cdot199\cdot419\cdot449\cdot521\cdot571
\cdot911\cdot2089\cdot2551\cdot5851\cdot11969$$
and
$$41\cdot67\cdot79\cdot181\cdot199\cdot233\cdot239\cdot307\cdot449\cdot521
\cdot1217\cdot1871\cdot4159\cdot5851\cdot9281.$$

We used a similar method to construct elements of $C(2,L_2)$.  We divided the set $P(2,L_2)$
into the set $S_1$ of its smallest $20$ members, the set $S_2$ of the $20$ smallest elements not
in $S_1$, and the set $S_3$ of the remaining $18$ elements,
and we defined $m_i$ as before.
We expect that there are about $2^{F(L_2)}\approx 2^{16.132}$ elements in $C(2,L_2)$,
so we expect that for every $4$ divisors $d$ of $m_3$ we should find one element
in $X\cap Y_d$.  This expectation is borne out by experimentation.  For example,
of the $18$ prime divisors of $m_3$, four give rise to Carmichael numbers of order $2$;
these Carmichael numbers are
$$23 \cdot 43 \cdot 59 \cdot 61 \cdot 79 \cdot 89 \cdot 113 \cdot 131 \cdot 151 \cdot 191 
\cdot 307 \cdot 311 \cdot 373 \cdot 419 \cdot 433 \cdot 463 \cdot 701 \cdot 1217 \cdot 2551,$$
$$23 \cdot 53 \cdot 59 \cdot 79 \cdot 89 \cdot 101 \cdot 109 \cdot 113 \cdot 131 \cdot 181 
\cdot 199 \cdot 233 \cdot 307 \cdot 349 \cdot 433 \cdot 701 \cdot 911 \cdot 1217 \cdot 4523,$$
\begin{multline*}
61 \cdot 67 \cdot 71 \cdot 89 \cdot 101 \cdot 103 \cdot 113 \cdot 151 \cdot 181 \cdot 191 \cdot 199 \cdot 233\\
\cdot 239 \cdot 271 \cdot 307 \cdot 419 \cdot 463 \cdot 521 \cdot 571 \cdot 701 \cdot 911 \cdot 5279,
\end{multline*}
and
\begin{multline*}
41 \cdot 43 \cdot 53 \cdot 61 \cdot 89 \cdot 103 \cdot 113 \cdot 151 \cdot 191 \cdot 311 \cdot 349\\
 \cdot 373 \cdot 419 \cdot 433 \cdot 463 \cdot 521 \cdot 571 \cdot 701 \cdot 929 \cdot 15313.
\end{multline*}

\section{Examples of non-rigid Carmichael numbers}
\label{S-examples}

Let $m$ be a positive integer.
Recall that we defined a rigid Carmichael number of order $m$ to be
a positive composite integer $n$ for which
$x\mapsto x^n$ is the identity map on every finite \'etale $\BZ/n\BZ$-algebra that
can be generated as a $\BZ/n\BZ$-module by $m$ elements.  Using arguments
like those in the proof of Theorem~\ref{characterization}, one can show that
a positive integer $n$ is a rigid Carmichael number of order $m$
if and only if $n$ is a squarefree composite integer such that $n\equiv 1\bmod (p^r-1)$
for every $r$ with $1\le r\le m$ and for every prime divisor $p$ of $n$.

We see that every element of the set $C(m,L)$ from Section~\ref{S-heuristics}
is a rigid Carmichael number of order $m$.
It is natural to ask whether all Carmichael numbers of
order $m$ are also rigid Carmichael numbers.  
The answer is no; we prove this by producing several Carmichael numbers $n$ of 
order $2$ each having a prime divisor $p$ with $n\not\equiv 1\bmod (p^2-1)$.

Let $L_0$ be a positive integer and let $p_0$ be a prime number that does not divide $L_0$ 
and such that $\gcd(L_0,p_0^2-1)$ divides $p_0 - 1$.  Let $P(2,L_0)$ be as in Section~\ref{S-heuristics},
and let $C(2,L_0,p_0)$ denote the set of integers of the form $p_0n_0$,
where $n_0$ is a squarefree integer, all of whose prime
factors lie in $P(2,L_0)$, such that $n_0 \equiv 1\bmod (p_0^2-1)$ and
$p_0n_0\equiv 1\bmod L_0$.  (Our assumption on $\gcd(L_0,p_0^2-1)$ ensures that
such $n_0$ are not barred from existence by congruence conditions.)  Then 
for every $n$ in $C(2,L_0,p_0)$ and every prime divisor $p$ of $n$ we have
$$ n \equiv 
\begin{cases}
1 \bmod (p^2-1) & \text{if $p\neq p_0$} \\
p \bmod (p^2-1) & \text{if $p=p_0$.}\\
\end{cases}$$
Since such an $n$ is squarefree, Theorem~\ref{characterization} shows that it is
a Carmichael number of order $2$, but it certainly is not a rigid Carmichael
number of order $2$.

If $L_0$ and $p_0$ are as above, let $L$ be the least common multiple
of $L_0$ and~$p_0^2 - 1$.  
Heuristics as in Section~\ref{S-heuristics} indicate that we should expect
there to be about $2^{\#P(2,L_0)}/\varphi(L)$ elements in the
set $C(2,L_0,p_0)$.  

For example, suppose we take $L_0$ to be
$2^7 \cdot 3^3 \cdot 5^2 \cdot 7 \cdot 11 \cdot 13 \cdot 17 \cdot 19 \cdot 29 \cdot 31$
(the number called $L_2$ in Section~\ref{S-constructing}), and suppose we let $p_0 = 1153$
(the smallest prime that does not divide $L_0$ and that satisfies the gcd condition mentioned
above).
Since $\#P(2,L_0) = 58$ and 
$\log \varphi(L) / \log 2 \approx 52$, we expect there to be about $64$ 
integers in $C(2,L_0,p_0)$.
We used a slightly modified version of the technique described in the preceding section
to search for elements of $C(2,L_0,p_0)$.  (We chose the subsets $S_1$ and $S_2$ of $P(2,L_0)$ so that
they each contained only quadratic residues modulo~$5$ ---
this allowed us to immediately disregard those divisors of $m_3$ that are quadratic residues 
modulo~$5$, since we were trying to find a divisor of $m_1m_2m_3$ that is congruent 
modulo $L$ to a quadratic nonresidue modulo~$5$.)
We found there to be $53$ elements in $C(2,L_0,p_0)$; the smallest of these is
\begin{multline*}
23\cdot 67\cdot 71\cdot 89\cdot 109\cdot 113\cdot 191\cdot 199\cdot 233\cdot 239\cdot 271\cdot 307\cdot 373\\
\cdot 419\cdot 521\cdot 911\cdot 929\cdot 1153\cdot 1217\cdot 1429\cdot 2089\cdot 2729\cdot 23561,
\end{multline*}
and the largest is
\begin{multline*}
23\cdot 37\cdot 43\cdot 53\cdot 59\cdot 61\cdot 67\cdot 71\cdot 89\cdot 103\cdot 109\cdot 113\cdot 131\cdot 181\cdot 191\cdot 199\cdot 239\cdot 271\\
\cdot 311\cdot 373\cdot 379\cdot 419\cdot 433\cdot 463\cdot 521\cdot 683\cdot 701\cdot 911\cdot 929\cdot 991\cdot 1153\cdot 1429\\
\cdot 2089\cdot 2551\cdot 3191\cdot 4159\cdot 5279\cdot 11969\cdot 15809\cdot 23561\cdot 23869\cdot 244529.
\end{multline*}

\end{document}